\documentclass[psamsfonts]{amsart}

\usepackage{amssymb,amsfonts,amscd, latexsym, graphicx, psfrag}
\usepackage[all,arc]{xy}
\usepackage{enumerate}
\usepackage{mathrsfs}
\newtheorem{thm}{Theorem}[section]

\newtheorem{prop}[thm]{Proposition}

\theoremstyle{definition}
\newtheorem{defn}[thm]{Definition}

\theoremstyle{remark}

\newtheorem{example}[thm]{Example}

\newcommand{\bC}{\mathbb{C}}

\newcommand{\bZ}{\mathbb{Z}}

\newcommand{\Ga}{\Gamma}

\newcommand{\GL}{\mathrm{GL}}

\newcommand{\rep}{\mathrm{Rep}}

\newcommand{\be}{\mathbf{e}}
\newcommand{\bu}{\mathbf{u}}

\newcommand{\cB}{\mathcal{B}}

\newcommand{\cD}{\mathcal{D}}
\newcommand{\cE}{\mathcal{E}}
\newcommand{\cF}{\mathcal{F}}
\newcommand{\cI}{\mathcal{I}}
\newcommand{\cL}{\mathcal{L}}
\newcommand{\cM}{\mathcal{M}}

\newcommand{\cR}{\mathcal{R}}
\newcommand{\cQ}{\mathcal{Q}}

\newcommand{\cX}{\mathcal{X}}

\newcommand{\sw}{\mathsf{w}}

\newcommand{\tcL}{\widetilde{\cL} }
\newcommand{\tT}{ \widetilde{\mathbb{T}} }

\newcommand{\hbu}{\hat{\bu}}

\newcommand{\age}{\mathrm{age}}

\newcommand{{\inv} }{\mathrm{inv}}
\newcommand{\ev}{\mathrm{ev}}
\newcommand{\Aut}{\mathrm{Aut}}
\newcommand{\rank}{\mathrm{rank}}
\newcommand{\val}{ {\mathrm{val}} }
\newcommand{\tw}{ {\mathrm{tw}} }
\newcommand{\Conj}{ \mathrm{Conj} }

\newcommand{\bGa}{\mathbf{\Ga}}

\newcommand{\Mbar}{\overline{\cM}}

\newcommand{\vGa}{\vec{\Ga}}

\newcommand{\C}{\mathbb{C}}

\newcommand{\Z}{\mathbb{Z}}

\newcommand{\BG}{\cB G}
\newcommand{\IBG}{\cI \cB G}
\newcommand{\IX}{\cI \cX}

\begin{document}
\title{Twisted Equivariant Gromov-Witten Theory of the Classifying Space
of a Finite Group}
\author{Zhuoming Lan}
\address{Zhuoming Lan, Department of Mathematical Sciences, Tsinghua University, Haidian District, Beijing 100084, China}
\email{lanzm21@mails.tsinghua.edu.cn}

\author{Zhengyu Zong}
\address{Zhengyu Zong, Department of Mathematical Sciences,
Tsinghua University, Haidian District, Beijing 100084, China}
\email{zyzong@mail.tsinghua.edu.cn}

\begin{abstract}
For any finite group $G$, the equivariant Gromov-Witten invariants of $[\bC^r/G]$ can be viewed as a certain twisted Gromov-Witten invariants of the classifying stack $\cB G$. In this paper, we
use Tseng's orbifold quantum Riemann-Roch theorem to express the equivariant Gromov-Witten invariants of $[\bC^r/G]$ as a sum over Feynman graphs, where the weight of each graph is expressed in terms of descendant integrals over moduli spaces of stable curves and representations of $G$.
    \end{abstract}

\maketitle

\tableofcontents

\section{Introduction}

\subsection{Historical background and motivation}
Gromov-Witten invariants are virtual counts of parametrized algebraic curves in algebraic manifolds, or more generally,
parametrized holomorphic curves in K\"{a}hler manifolds. Gromov-Witten theory can be viewed
as a mathematical theory of A-model topological string theory. The topological string theory on orbifolds was constructed
decades ago by physicists \cite{DHVW85, DHVW86}, and many works followed in both mathematics and physics (e.g. \cite{HV, DFMS, BR, Kn, Roan, CV}).  Orbifold quantum cohomology was discussed in \cite{Za} along with many examples, in both abelian and non-abelian quotients. Later, the mathematical definition of orbifold Gromov-Witten theory and quantum cohomology was laid in \cite{CR02, AGV02, AGV08}.

Orbifold Gromov-Witten theory appears naturally in enumerative geometry and mathematical physics. One interesting case is the Remodeling Conjecture proposed by Bouchard-Klemm-Mari\~{n}o-Vafa \cite{BKMP09, BKMP10}. The Remodeling Conjecture is a version of all genus open-closed mirror symmetry for semi-projective toric Calabi-
Yau 3-orbifolds. The B-model is given by Eynard-Orantin topological recursion \cite{EO07} on the mirror curve. The open string sector of the Remodeling Conjecture for $\bC^3$ was proved independently
by L. Chen \cite{Ch09} and J. Zhou \cite{Zh09}; the closed string sector of the Remodeling
Conjecture for $\bC^3$ was proved independently by Bouchard-Catuneanu-
Marchal-Sulkowski \cite{BCMS} and S. Zhu \cite{Zhu}. Eynard and Orantin provided a proof of
the Remodeling Conjecture for smooth semi-projective toric Calabi-Yau 3-folds in \cite{EO12}. The Remodeling Conjecture for general semi-projective toric Calabi-
Yau 3-orbifolds is proved in \cite{FLZ2}.

In \cite{FLZ1}, the Remodeling Conjecture is proved for affine toric Calabi-Yau 3-
folds. The result in \cite{FLZ1} plays a crucial role in the proof of the Remodeling Conjecture for general semi-projective toric Calabi-Yau 3-orbifolds given in \cite{FLZ2}. On the other hand, a key ingredient in the proof of \cite{FLZ1} is the study of equivariant Gromov-Witten theory of affine smooth toric Deligne-Mumford stacks given in \cite{FLZ0}. An affine smooth toric Deligne-Mumford stack can be written as the quotient stack $[\bC^r/G]$, where $G$ is a finite abelian group. The equivariant Gromov-Witten theory of $[\bC^r/G]$ can be viewed as a certain twisted Gromov-Witten theory \cite{Ts10, CG07} of the classifying space $\cB G$. By Tseng's orbifold quantum Riemann-Roch theorem \cite{Ts10}, the equivariant Gromov-Witten theory of $[\bC^r/G]$ can be expressed as a sum over stable Feynman graphs, where the weights assigned to vertices of these graphs are essentially descendant integrals over $\Mbar_{g,n}$. When $r=3$, this graph sum formula can be applied to the proof of the Remodeling Conjecture for affine toric Calabi-Yau 3-folds \cite{FLZ1} by comparing with the B-model graph sum formula for Eynard-Orantin invariants \cite{DOSS}.

In \cite{bb}, the all genus mirror symmetry for $[\bC^3/G]$ is conjecturally generalized to the case where $G$ is possibly non-abelian. A natural question is that whether one can generalize the proof in \cite{FLZ1} to the non-abelian case. The aim of this paper is to study the equivariant Gromov-Witten theory of the quotient stack $[\bC^r/G]$, where $G$ is an arbitrary finite group, and obtain the corresponding graph sum formula. The result in this paper can be viewed as a non-abelian generalization of \cite{FLZ0} and the authors hope that this graph sum formula can be applied to give a proof of the all genus mirror symmetry conjectured in \cite{bb}.

\subsection{Statement of the main result}
Let $G$ be a finite group and we consider an action $G$ on $\bC^r$ via a representation $\rho$ of $G$ on $\bC^r$. We decompose $\rho=\sum_{i=1}^{m}\rho_{i}$ into irreducible representations, with $\dim \rho_{i}=r_i$, and let $\C^{r}=\bigoplus_{i=1}^{m}\C^{r_i}$ be the corresponding decomposition of $\bC^r$ into $G$ invariant subspaces. Consider the action of $\tT=(\C^{*})^{m}$ on $[\C^{r}/G]$, induced by the following action on $\C^{r}=\bigoplus_{i=1}^{m}\C^{r_i}$:
$$\tT\times\bigoplus_{i=1}^{m}\C^{r_i}\to\bigoplus_{i=1}^{m}\C^{r_i}$$
$$
((z_{1},\dots,z_{m}),(\mathbf{v}_{1},\dots,\mathbf{v}_{m}))\mapsto(z_{1}\mathbf{v}_{1},\dots,z_{m}\mathbf{v}_{m}).
$$

We will use Tseng's orbifold quantum Riemann-Roch theorem \cite{Ts10} to obtain the following main theorem (See Theorem \ref{main}).

\begin{thm}
The $\tT$ equivariant Gromov-Witten theory of $[\bC^r/G]$ can be expressed as a sum over stable Feynman graphs, where the weights assigned to vertices of these graphs are essentially descendant integrals over $\Mbar_{g,n}$.
\end{thm}

\subsection{Future work}
In \cite{bb}, an all genus mirror symmetry of $[\bC^3/G]$ is conjectured for some non-abelian groups $G$'s. The B-model is conjecturally given by the Eynard-Orantin topological recursion on the mirror curve. One can apply the graph sum formula in this paper to the above cases and compare it with the B-model graph sum formula for Eynard-Orantin invariants \cite{DOSS}. We expect that one would be able to match these two graph sum formulas and therefore prove the conjecture in \cite{bb}.

\subsection{Overview of the paper}
In Section 2, we study the Chen-Ruan orbifold cohomology of $\BG$ and $[\C^r/G]$. In Section 3, we study the equivariant Gromov-Witten theory of the quotient stack $[\bC^r/G]$ and explain it as a certain twisted Gromov-Witten theory of $\cB G$. In Section 4, we use Tseng's orbifold quantum Riemann-Roch theorem \cite{Ts10} to obtain the graph sum formula for the equivariant Gromov-Witten theory of $[\bC^r/G]$.

\subsection*{Acknowledgement}
The authors wish to thank Chiu-Chu Melissa Liu and Bohan Fang for useful discussions. The work of the second author is partially supported by NSFC grant No. 11701315.

\section{Geometry and Chen-Ruan orbifold cohomology of $[\C^r/G]$}
\subsection{Chen-Ruan orbifold cohomology of $\BG$}
Let $G$ be a finite group and $\BG$ the classifying stack of $G$. Let $\Conj(G)$ be the set of conjugacy classes of $G$. The inertia stack of $\BG$ is
$$
\IBG = \bigcup_{[h]\in \Conj(G)}(\BG)_h,
$$
and
$$
(\BG)_h =[\{h\}/C(h)]\cong \BG
$$
where $C(h)$ is the centralizer of $h$ in $G$.
As a graded vector space over $\C$, the Chen-Ruan orbifold
cohomology of $\BG$ is
$$
H^*_\mathrm{CR}(\BG;\C) = H^*(\IBG;\C) = \bigoplus_{[h]\in G} H^0((\BG)_h;\C),
$$
where $H^0(\BG_h;\C)=\C 1_h$.
The orbifold Poincar\'{e} pairing of $H^*_\mathrm{CR}(\BG;\C)$ is given by
$$
\langle 1_h,
 1_{h'}\rangle = \frac{\delta_{h^{-1},h'}}{|C(h)|}.
$$
The orbifold cup product of $H^*_\mathrm{CR}(\BG;\C)$ is given by
$$
1_h \star 1_{h'} =\sum_{g\in[h],g'\in[h']}\frac{|C(gg')|1_{gg'}}{|G|}.
$$

Following \cite{JK}, we define a canonical basis for
the semisimple algebra $H^*_\mathrm{CR}(\BG;\C)$.
Define $\rep(G)$ to be the set of irreducible representations. Given $\gamma\in \rep(G)$, define
$$
\phi_\gamma :=\frac{|V_{\gamma}|}{|G|}\sum_{[h]\in \Conj(G)}\chi_\gamma(h^{-1}) 1_h
$$
where $V_{\gamma}$ is the vector space corresponding to $\gamma$ and $|V_{\gamma}|$ is the dimension of $V_\gamma$.

Then
$$
H^*_\mathrm{CR}(\BG;\C) =\bigoplus_{[h]\in \Conj(G)}\C 1_h  =\bigoplus_{\gamma \in \rep(G)}\C \phi_\gamma.
$$
Recall that we have the orthogonality of characters:
\begin{enumerate}
\item For any $\gamma,\gamma'\in \rep(G)$,
$\displaystyle{
\frac{1}{|G|} \sum_{h\in G} \chi_{\gamma}(h^{-1})\chi_{\gamma'}(h)=\delta_{\gamma,\gamma'}
},$
\item For any $h, h'\in G$,
$\displaystyle{
\frac{1}{|G|}\sum_{\gamma\in \rep(G)}\chi_\gamma(h^{-1})\chi_\gamma(h') =\delta_{[h],[h']}.
}$
\end{enumerate}
Therefore, by defining $\nu_{\gamma}=(\frac{\dim V_{\gamma}}{|G|})^{2}$, we have
$$
\langle \phi_{\gamma}, \phi_{\gamma'}\rangle = \nu_{\gamma}\delta_{\gamma,\gamma'},
$$
and
$$
\phi_{\gamma}\star\phi_{\gamma'}= \delta_{\gamma, \gamma'}\phi_{\gamma}.
$$
This shows that $\{\phi_\gamma: \gamma\in \rep(G)\}$ is a canonical basis
of the semisimple $\C$-algebra $H^*_\mathrm{CR}(\cB G;\C)$.

\subsection{Equivariant Chen-Ruan orbifold cohomology of $\cX=[\C^r/G]$} \label{eqbg}
Let $G$ be a finite group and we consider an action $G$ on $\bC^r$ via a representation $\rho$ of $G$ on $\bC^r$. We decompose $\rho=\sum_{i=1}^{m}\rho_{i}$ into irreducible representations, with $\dim \rho_{i}=r_i$, and let $\C^{r}=\bigoplus_{i=1}^{m}\C^{r_i}$ be the corresponding decomposition into $G$ invariant subspaces. Given any $h\in G$ and $i\in \{1,\ldots, r\}$, consider $\rho(h)\in \bigoplus_{\rho_{i}}\GL_{r_{i}}(\C)$ and we can diagonalize it to be
$$\bigoplus_{i=1}^{m}\mathbf{diag}(\exp(2\pi ih_{i1}),\dots,\exp(2\pi ih_{ir_{i}}))$$
where $h_{ij}\in[0,1)$.

Define the age of $h$ with respect to $\rho_{i}$ to be
$$
\age_{\rho_i}(h)=\sum_{j=1}^{r_i}h_{ij}
$$
and the age with respect to $\rho$ to be

$$
\age_{\rho}(h)=\sum_{i=1}^{m}\age_{\rho_{i}}(h)
$$

Let $(\C^r)^h$ denote the $h$-invariant subspace of $\C^r$. Then
$$
\dim_\C (\C^r)^h  = \sum_{i=1}^{m}\sum_{j=1}^{r_{i}} \delta_{h_{ij},0}.
$$
The inertial stack of $\cX$ is
$$
\IX  = \bigcup_{[h]\in G}\cX_h,
$$
where
$$
\cX_h =[(\C^r)^h/C(h)].
$$
In particular,
$$
\cX_1 =[\C^r/G]=\cX.
$$
As a graded vector space over $\C$,
$$
H^*_\mathrm{CR}(\cX;\C) = \bigoplus_{h\in G} H^*(\cX_h; \C)[2\age_{\rho}(h)] =\bigoplus_{h\in G} \C 1_h,
$$
where $\deg(1_h)=2\age_{\rho}(h)\in \mathbb{Q}$. The orbifold Poincar\'{e} pairing of the (non-equivariant) Chen-Ruan orbifold cohomology $H_\mathrm{CR}^*(\cX;\C)$ is given by
$$
\langle 1_{h}, 1_{h'}\rangle_{\cX} =\frac{1}{|C(h)|}\delta_{h^{-1},h'}
\cdot \delta_{0,\dim_\C \cX_h}
$$

Consider the action of $\tT=(\C^{*})^{m}$ on $[\C^{r}/G]$, induced by the following action on $\C^{r}=\bigoplus_{i=1}^{m}\C^{r_i}$:
$$\tT\times\bigoplus_{i=1}^{m}\C^{r_i}\to\bigoplus_{i=1}^{m}\C^{r_i}$$
$$
((z_{1},\dots,z_{m}),(\mathbf{v}_{1},\dots,\mathbf{v}_{m}))\mapsto(z_{1}\mathbf{v}_{1},\dots,z_{m}\mathbf{v}_{m}).
$$
Let $\cR=H^*(B\tT;\C)=\C[\sw_{\rho_1},\dots,\sw_{\rho_m}]$, where
$\sw_{\rho_1},\dots,\sw_{\rho_m}$ are the first Chern classes of the universal line bundles over $\cB\tT$.
The $\tT$-equivariant Chen-Ruan orbifold cohomology $H^*_{\mathrm{CR},\tT}([\C^r/G];\C)$
is an $\cR$-module.  Given $h\in G$, define
$$
\mathbf{e}_h := \prod_{i=1}^m\prod_{j=1}^{r_{i}} \sw_{\rho_{i}}^{\delta_{h_{ij},0}} \in \cR.
$$
In particular,
$$
\mathbf{e}_1 =\prod_{i=1}^{m}\sw_{\rho_{i}}^{\dim\rho_{i}}.
$$
Then the $\tT$-equivariant Euler class of $0_h:=[0/G]$ in $\cX_h = [(\C^r)^h/G]$ is
$$
e_{\tT}(T_{0_h}\cX_h) = \mathbf{e}_h 1_h \in H^*_{\tT}(\cX_h;\C) = \cR 1_h.
$$

Define
$$
\cR' =\C[\sw_{\rho_{1}}^{1/|G|},\dots,\sw_{\rho_{m}}^{1/|G|}],
$$
which is a finite extension of $\cR$. Let $\cQ$ and $\cQ'$ be the fractional fields of
$\cR$ and $\cR'$, respectively. The $\tT$-equivariant Poincar\'{e} pairing of
$H^*_{\tT,\mathrm{CR}}(\cX;\C)\otimes_{\cR}\cQ$ (which is isomorphic to  $H^*_{\mathrm{CR}}(\cX;\cQ)$
as a vector space over $\cQ$) is given by
$$
\langle 1_{h}, 1_{h'}\rangle_{\cX} =\frac{1}{|G|}\cdot\frac{\delta_{h^{-1},h'} }{\mathbf{e}_h}\in \cQ.
$$
The $\tT$-equivariant orbifold cup product  of $H^*_{\tT,\mathrm{CR}}(\cX;\C)\otimes_{\cR}\cQ$ is given by
$$
1_h \star_{\cX} 1_{h'} = \sum_{g\in[h],g'\in[h']}\frac{|C(gg')|\prod_{i=1}^{m}\sw_{\rho_{i}}^{\age_{\rho_{i}}(g)+\age_{\rho_{i}}(g')-\age_{\rho_{i}}(gg')}}{|G|}1_{gg'}.
$$

Define
\begin{equation}\label{eqn:bar-one}
\bar{1}_h:= \frac{1_h}{\prod_{i=1}^{m}\sw_{\rho_{i}}^{age_{\rho_{i}}(h)}} \in H^*_{\tT,\mathrm{CR}}(\cX;\C)\otimes_{\cR}\cQ'.
\end{equation}
Then
$$
\langle \bar{1}_h, \bar{1}_{h'}\rangle_{\cX} =
\frac{\delta_{h^{-1},h'}}{|G|\mathbf{e}_1}
$$
and
$$
\bar{1}_h \star_{\cX}\bar{1}_{h'} =\sum_{g\in[h],g'\in[h']}\frac{|C(gg')|}{|G|}\bar{1}_{gg'}.
$$

Given $\gamma\in \rep(G)$, define
\begin{equation}\label{eqn:bar-phi}
\bar{\phi}_\gamma := \frac{|V_{\gamma}|}{|G|}\sum_{[h]\in G} \chi_\gamma(h^{-1}) \bar{1}_h.
\end{equation}

Then by defining $\bar{\nu_{\gamma}}=(\frac{|V_{\gamma}|}{G\sqrt{\mathbf{e}_{1}}})^{2}$
$$
\langle \bar{\phi}_{\gamma},\bar{\phi}_{\gamma'}\rangle_{\cX} = \bar{\nu}_{\gamma}\delta_{\gamma\gamma'}
$$
and
$$
\bar{\phi}_{\gamma}\star_{\cX} \bar{\phi}_{\gamma'} =  \delta_{\gamma \gamma'} \bar{\phi}_{\gamma}.
$$

In particular, $\{\bar{\phi}_\gamma:\gamma\in \rep(G)\}$ is a canonical basis
of the semisimple $\cQ'$-algebra $H^*_{\tT,\mathrm{CR}}(\cX;\C)\otimes_{\cR}\cQ'$.

\section{Equivariant Gromov-Witten theory of $[\C^{r}/G]$}

\subsection{The total descendant Gromov-Witten potential of $\BG$}
In \cite{JK}, Jarvis-Kimura studied all genus
descendant orbifold Gromov-Witten invariants of $\BG$
for any finite group $G$. In this subsection, we recall
their results.

Let $\Mbar_{g,n}(\BG)$ be moduli space of genus $g$, $n$-pointed
twisted stable maps to $\BG$, and let
$$
\ev_j:\Mbar_{g,n}(\BG)\to \IBG = \bigcup_{[h]\in G}(\BG)_h
$$
be the evaluation at the $j$-th marked point. Given
$h_1,\ldots, h_n\in G$, define
$$
\Mbar_{g,(h_1,\ldots, h_n)}(\BG):= \bigcap_{j=1}^n \ev_j^{-1}\left((\BG)_{h_j}\right).
$$
Then $\Mbar_{g,(h_1,\ldots, h_n)}(\BG) $ is empty unless there exists $g_{i}\in[h_{i}],i=1,\dots n$ such that $g_1\cdots g_n=1$.

Given $\beta_1,\ldots,\beta_n\in H^*_\mathrm{CR}(\BG;\C)$
and $a_1,\ldots, a_n\in \Z_{\geq 0}$, define the correlator
$$
\langle \tau_{a_1}(\beta_1) \cdots \tau_{a_n}(\beta_n)\rangle_{g,n}^{\BG}:=
\int_{\Mbar_{g,n}(\BG)}\prod_{j=1}^n \bar{\psi}_j^{a_j}\ev_j^*\beta_j.
$$
With this notation, we have
\begin{itemize}
\item For any $h_1,\ldots, h_n\in G$,
\begin{eqnarray*}
&& \langle \tau_{a_1}(1_{h_1})\cdots \tau_{a_n}(1_{h_n})\rangle_{g,n}^{\BG}\\
&=& \begin{cases}
\Omega^{G}_{g,n}(\gamma)
 \int_{\Mbar_{g,n}}\psi_1^{a_1}\cdots \psi_n^{a_n}&
\textup{if }\Mbar_{g,(h_1,\ldots, h_n)}(\BG) \textup{ is nonempty}\textup{ and } \sum_{i=1}^n a_i =3g-3+n,\\
0 &\textup{otherwise}.
\end{cases}
\end{eqnarray*}
where $$\Omega^{G}_{g,n}(\gamma)=|X^{G}_{g,n}(\gamma)|/|G|$$ and $$X^{G}_{g,n}(\gamma)= \{(\alpha_1,\dots, \alpha_g, \beta_1,\dots,\beta_g, \sigma_1,\dots, \sigma_n)|\prod^{g}_{i=1}[\alpha_i,\beta_{i}] = \prod^{n}_{j=1}\sigma_{j}, \sigma_{j}\in [\gamma_{j}]\}$$
\item For any $\gamma_1,\ldots, \gamma_n\in \rep(G)$,
\begin{eqnarray*}
&& \langle \tau_{a_1}(\phi_{\gamma_1})\cdots \tau_{a_n}(\phi_{\gamma_n}) \rangle_{g,n}^{\BG} \\
&=& \begin{cases}
\nu_{\gamma}^{1-g}\int_{\Mbar_{g,n}}\psi_1^{a_1}\cdots \psi_n^{a_n}&
\textup{if }\gamma_1 =\cdots =\gamma_n \textup{ and } \sum_{i=1}^n a_i =3g-3+n,\\
0 &\textup{otherwise}.
\end{cases}
\end{eqnarray*}
\end{itemize}

We have the following two propositions.
\begin{prop}[\cite{JK}]
The correlators $<\tau_{a_{1}}(e[[\sigma_1]])\dots  \tau_{a_n}
(e[[\sigma_n]])>^{\BG}_{g,n}$ are related to the
usual correlators  $<\tau_{a_{1}}\dots  \tau_{a_n}>_{g}$ (corresponding to the case of $G = \{1\}$) by
$<\tau_{a_{1}}(e[[\sigma_1]])\dots  \tau_{a_n}
(e[[\sigma_n]])>^{\BG}_{g,n}=<\tau_{a_{1}}\dots  \tau_{a_n}>_{g,n} \Omega^{G}_{g,n}(\gamma)$
where $\Omega^{G}_{g,n}(\gamma)=|X^{G}_{g,n}|(\gamma)/|G|$ and $X^{G}_{g,n}(\gamma)= \{(\alpha_1,\dots, \alpha_g, \beta_1,\dots,\beta_g, \sigma_1,\dots, \sigma_n)|\prod^{g}_{i=1}[\alpha_i,\beta_{i}] = \prod^{n}_{j=1}\sigma_{j}, \sigma_{j}\in [\gamma_{j}]\}.$

\end{prop}

The formula of correlation function for $\BG$ under Frobenius basis is:

\begin{prop}[\cite{JK}]
 We have   $<\tau_{a_1}(\phi_{\alpha})\dots  \tau_{a_n}(\phi_{\alpha})>_{g,n}^{BG}=\nu_{\alpha}^{1-g}<\tau_{a_1}\dots  \tau_{a_{n}}>_{g,n}$ and other types of correlators are equal to 0.
     Here $\nu_{\alpha}=(\frac{|V_{\alpha}|}{|G|})^{2}$.
\end{prop}

\subsection{The total descendant equivariant Gromov-Witten potential of $\mathcal{X}=[\C^{r}/G]$}
Given a representation $\rho$ of $G$, we denote the vector bundle over $\BG$ corresponding to $\rho$ as $F_{\rho}$. Then the quotient stack $\mathcal{X}=[\C^{r}/G]$ is the total space of the vector bundle of the form
$$F_{\rho}=\bigoplus_{i=1}^{m}F_{\rho_i},\rho_{i}\in\rep(G).$$

The vector bundle $F_{\rho}\to \BG$ is equipped with a $\widetilde{\mathbb{T}}=(\C^{*})^{m}$-equivariant structure as defined in \ref{eqbg}. Let
$$F_{g,n}^{\rho}=R^{*}\pi_{*}f^{*}F_{\rho}=R^{0}\pi_{*}f^{*}F_{\rho}-R^{1}\pi_{*}f^{*}F_{\rho}\in K_{\widetilde{\mathbb{T}}}(\overline{\mathcal{M}}_{g,n}(BG))$$
be the K-theoretic pushforward, where $K_{\widetilde{\mathbb{T}}}$ denotes the $\widetilde{\mathbb{T}}$-equivariant K-theory,
$$\pi: C_{g,n}\to\overline{\mathcal{M}}_{g,n}(\BG)
$$
is the projection from the universal curve, and
$$f:C_{g,n}\to \BG$$ is the universal map.

By definition, $F^{\rho}_{g,n}$ has the natural direct sum decomposition:

$$F^{\rho}_{g,n}=\bigoplus_{i=1}^{m}F^{\rho_i}_{g,n}.$$

By orbifold Riemann-Roch, we have:
\begin{equation} \label{eqn:rankF}
 \mathrm{rank} F_{g,n}^{\rho_i}|_{\overline{\mathcal{M}}_{g,n}(h_{1},\dots,h_{n})}=(1-g)\dim \rho_{i}-\sum_{j=1}^{n}age_{\rho_i}(h_{i})
\end{equation}
\begin{equation}
\mathrm{rank} F_{g,n}^{\rho}|_{\overline{\mathcal{M}}_{g,n}(h_{1},\dots,h_{n})}=(1-g)\dim \rho-\sum_{j=1}^{n}age_{\rho}(h_{i}).
\end{equation}

The $\tT$-equivariant Euler class of $F^{\rho}_{g,n}$ is
\begin{equation}\label{eqn:eT}
e_{\tT}(F_{g,n}^{\rho}) =\prod_{i=1}^m \left(\sw_{\rho_i}^{\rank F^{\rho_i}_{g,n}}c_{\frac{1}{\sw_{\rho_i}}}(F^{\rho_i}_{g,n})\right),
\end{equation}
where
$$
c_t(F) = 1 + c_1(F) t + c_2(F) t^2 + \cdots
$$

Define the $\widetilde{\mathbb{T}}$-equivariant genus $g$ descendant orbifold Gromov-Witten invariants of $\mathcal{X}$ by
\begin{equation}\label{eqn:X}
\langle \tau_{a_1}(\beta_1) \cdots \tau_{a_n}(\beta_n)\rangle^{\cX}_{g,n}:=
\int_{\Mbar_{g,n}(\BG)}\frac{\prod_{j=1}^n \bar{\psi}_j^{a_j}\ev_j^*\beta_j}{e_{\tT}(F_{g,n}^{\rho})}.
\end{equation}

Here we use the $\cQ'$-vector space isomorphism
$H^*_{\tT,\mathrm{CR}}(\cX;\C)\otimes_{\cR}\cQ' \cong H^*_\mathrm{CR}(\BG;\cQ')$ to view
$\beta_1,\ldots,\beta_n$ as elements in $H^*_\mathrm{CR}(\BG;\cQ')$.
Introducing formal variables
$$
\hat{\bu}=\sum_{a\geq 0}\hat u_a z^a,
$$
where $\hat u_a \in H^*_\mathrm{CR}(\cX;\cQ)$, and define
$$
\langle \hat\bu,\ldots, \hat\bu \rangle^{\cX}_{g,n}
:= \sum_{a_1,\ldots, a_n\in \Z_{\geq 0}}  \langle\tau_{a_1}(\hat u_{a_1})\cdots, \tau_{a_n}(\hat u_{a_n}) \rangle^{\cX}_{g,n}.
$$
Define the generating functions
$$
\cF_g^{\cX}(\hat\bu):= \sum_{n\geq 0} \frac{1}{n!} \langle \hat\bu,\ldots, \hat\bu\rangle^{\cX}_{g,n},
$$
$$
\cD^{\cX}(\hat\bu,\hbar): =\exp\biggl(\sum_{g\geq 0}\hbar^{g-1}\cF_g^{\cX}(\hat\bu)\biggr).
$$
The generating function $\cD^{\cX}(\hat\bu,\hbar)$ is called the total descendant
$\tT$-equivariant Gromov-Witten potential of $\cX$.

\subsection{The twisted total descendant Gromov-Witten potential}
In this subsection, we define a twisted total descendant
Gromov-Witten potential $\cD^{\tw}(\bu, \hbar)$, which is related to $D^{\cX}(\hat\bu,\hbar)$
by a simple change of variables. In Section \ref{sec:graph}, we will write
$\cD^{\tw}(\bu,\hbar)$ as a graph sum.

Given $\beta_1,\ldots, \beta_n \in H^*_\mathrm{CR}(\BG;\cQ')$
and $a_1,\ldots, a_n \in \Z_{\geq 0}$, define twisted correlators:
\begin{equation}\label{eqn:twisted}
\langle \tau_{a_1}(\beta_1) \cdots \tau_{a_n}(\beta_n)\rangle^\tw_{g,n}:=
\int_{\Mbar_{g,n}(\BG)}\frac{\prod_{j=1}^n \bar{\psi}_j^{a_j}\ev_i^*\beta_j}{
\prod_{i=1}^r c_{\frac{1}{\sw_{\rho_i}}}(F^{\rho_{i}}_{g,n})}.
\end{equation}
By \eqref{eqn:rankF}, \eqref{eqn:eT}, \eqref{eqn:X}, and \eqref{eqn:twisted},
for any $h_1,\ldots, h_n\in G$,
$$
\langle \tau_{a_1}(1_{h_1}) \cdots \tau_{a_n}(1_{h_n})\rangle^\cX_{g,n}
=\prod_{i=1}^r \sw_{\rho_i}^{\dim \rho_i (g-1)+\sum_{j=1}^n \age_{\rho_i}(h_j)} \langle \tau_{a_1}(1_{h_1}) \cdots \tau_{a_n}(1_{h_n})\rangle^\tw_{g,n}
$$
or equivalently,
\begin{equation}
\langle \tau_{a_1}(\bar{1}_{h_1}) \cdots \tau_{a_n}(\bar{1}_{h_n})\rangle^\cX_{g,n}
=\be_1^{g-1} \langle \tau_{a_1}(1_1) \cdots \tau_{a_n}(1_n)\rangle^\tw_{g,n},
\end{equation}
where $\bar{1}_h$ is defined by \eqref{eqn:bar-one}.

Let $\phi_\gamma$ be defined as in \eqref{eqn:bar-phi}. Then
for any $\gamma_1,\ldots, \gamma_n\in \rep(G)$,
\begin{equation}
\langle \tau_{a_1}(\bar{\phi}_{\gamma_1}) \cdots \tau_{a_n}(\bar{\phi}_{\gamma_n})\rangle^\cX_{g,n}
=\be_1^{g-1} \langle \tau_{a_1}(\phi_{\gamma_1}) \cdots \tau_{a_n}(\phi_{\gamma_n})\rangle^\tw_{g,n}.
\end{equation}

Introduce formal variables
$$
\bu=\sum_{a\geq 0}u_a z^a,
$$
where $u_a \in H^*_\mathrm{CR}(\cX;\cQ)$, and define
$$
\langle \bu,\ldots \bu\rangle^\tw_{g,n}
=\sum_{a_1,\ldots, a_n\in \Z_{\geq 0}}
\langle \tau_a(u_{a_1}) \cdots \tau_{a_n}(u_{a_n})\rangle^{\tw}_{g,n}.
$$
Define the generating functions
$$
\cF_g^\tw(\bu):= \sum_{n\geq 0} \frac{1}{n!} \langle \bu,\ldots, \bu\rangle^\tw_{g,n}
$$
$$
\cD^\tw(\bu,\hbar): =\exp\biggl(\sum_{g\geq 0}\hbar^{g-1}\cF_g^\tw(\bu)\biggr).
$$
Then
$$
\cD^\cX(\hat\bu= \sum_{a\geq 0} \hat{u}_a z^a ,\hbar)
= \cD^\tw(\bu= \sum_{a\geq 0} u_a z^a,  \be_1 \hbar),
$$
where
$$
\hat{u}_a = \sum_{\alpha\in \rep(G)} u_a^\alpha\bar{\phi}_\alpha,
\quad u_a = \sum_{\alpha\in \rep(G)} u_a^\alpha \phi_\alpha.
$$
\section{Orbifold quantum Riemann-Roch theorem and graph sum formula} \label{sec:graph}
\subsection{Orbifold quantum Riemann-Roch theorem for $[\C^{r}/G]$}

In this section we apply the orbifold quantum Riemann-Roch in \cite{Ts10} to $[\C^{r}/G]$ case, which gives the relation of twisted descendant potentials with descendant potentials.

Let $\mathbf{c}(\mathbf{\cdot})$ be a characteristic class of $\BG$ and $F$ be a vector bundle over $\BG$. Define the twisted potential function to be:
\begin{equation}
    \cD_{(\mathbf{c}(\mathbf{\cdot}),F)}=\exp\biggl(\sum_{g\geq 0}\hbar^{g-1}\cF_g^{(\mathbf{c}(\mathbf{\cdot}),F)}(\bu)\biggr)
\end{equation}
where $\cF_g^{(\mathbf{c}(\mathbf{\cdot}),F)}(\bu)$ is defined similarly as $\cF_g^{\tw}(\bu)$ by taking the correlation function to be
\begin{equation}
    \langle \tau_{a_1}(\beta_1) \cdots \tau_{a_n}(\beta_n)\rangle^{(\mathbf{c}(\mathbf{\cdot}),F)}_{g,n}:=
\int_{\Mbar_{g,n}(\BG)}\frac{\prod_{j=1}^n \bar{\psi}_j^{a_j}\ev_i^*\beta_j}{
\mathbf{c}(R^{*}\pi_{*}f^{*}F)}.
\end{equation}

\begin{prop} [\cite{Ts10}] (Orbifold Quantum Riemann-Roch)
  Let $\cX=\BG$. Let $\cD_{s}=\cD_{(\mathbf{c}=\sum_{k}s_{k}ch_{k}(\mathbf{\cdot}),F)}$ with $s_{0}=0$. Let $A^{\wedge}$ be the Givental quantization defined in \cite{Gi3} of operator $A$.
Then we have
\begin{equation} \label{OQRR}
\begin{split}
&\quad \cD_s=\exp\left(\sum_{k\geq 1} s_k\left(\sum_{m>0} \frac{(A_m)_{k+1-m}z^{m-1}}{m!}\right)^\wedge\right)\exp\left(\sum_{k\geq 0} s_k\left(\frac{(A_0)_{k+1}}{z}\right)^\wedge\right)\cD^\cX
\end{split}
\end{equation}
where $(A_{m})_{k}=\sum_{\lambda}ch_{k}(F^{\lambda})B_{m}(\lambda)$ and $B_{m}(\lambda)$ is the Bernoulli polynomial. Here $F^\lambda$ is the sub-bundle of $F$ with $G$-eigenvalue $\exp(2\pi i\lambda)$ where $\lambda\in[0,1)$.
\end{prop}

The bundle $F$ in our case is $F_{\rho}$. The Chern character is as following:

$ch_{k}(F^{\lambda})=0$, for $k\neq 0$.

$ch_{k}(F^{\lambda})=\dim F^\lambda$, for $k=0$.

We decompose the bundle $F_{\rho}=\bigoplus F_{\rho_{i}}$. Define
$$(A_{t}^{i})_{k}=\sum_{\lambda}ch_{k}(F_{\rho_i}^{\lambda})B_{t}(\lambda).$$
Then we split $A_{t}$ into  the sum of $A_{t}^{i}$ by:
\begin{equation}
    A_{t}=\sum_{i=1}^{m}A_{t}^{i}
\end{equation}

The explicit expression of $A_{t}^{i}$ is as following.
\begin{defn}
    Let $F=F_{\rho}$. $A_{t}^{i}$ is the mapping:

    $1_{h}\to \sum_{l=0}^{o(h)-1} B_{t}(l/o(h))D_{\rho_{i}}^{h}(l)1_{h}$

    where $o(h)$ is the order of $h$ and $D_{\rho_i}^{h}(l)$ is the dimension of the bundle $F_{\rho_i}^{l/o(h)}$.
\end{defn}

   For simplicity, we define:

$$B^{\rho_i}_{t}(h)=\sum_{l=0}^{o(h)-1}B_{t}(l/o(h))D_{\rho_i}^{h}(l).$$

Then we calculate $s_{k}$ in \eqref{OQRR}.

Using the fact that $\exp(\sum_{k=1}^{\infty}(k-1)!(-1)^{k-1}ch_{k}(V)t^{k})=c_{t}(V)$ for any vector bundle $V$, for $s_{k}^{i}=(k-1)!(-w_{\rho_{i}})^{-k}$, $s_{0}=0$
we have:

 $\exp(\sum_{i=1}^{m}\sum_{k=1}^{\infty}s_{k}^{i}ch_{k}(F_{g,n}^{\rho_{i}}))=\frac{1}{\prod_{i=1}^{m}c_{\frac{1}{w_{\rho_{i}}}}(F^{\rho_{i}}_{g,n})}.$

Now we get the quantization formula for $D^{tw}$ by setting $A_{t}^{i}$ and $s_{k}^{i}$ as above.
\begin{equation}\label{eqn:QRR}
\cD^{\tw} = \exp\Bigl(\sum_{i=1}^{m}\sum_{t=1}^{\infty}\frac{(-1)^{t}}{t(t+1)}A^{i}_{t+1}(\frac{z}{\sw_{\rho_{i}}})^{t}\Bigr)^{\wedge}
\cD^{BG}.
\end{equation}

Under the canonical basis $\{\phi_{\alpha}\}$, the operator $A_{t}^{i}$ is expressed in the following form:
$$A_{t}^{i}(\phi_{\beta})=\sum_{\alpha}(E^{i}_{t})^{\alpha}_{\beta}\phi_{\alpha},$$
where
$$(E^{i}_{t})^{\alpha}_{\beta}=\frac{|V_{\alpha}|}{|G||V_{\beta}|}\sum_{h}\chi_{\alpha}(h)\chi_{\beta}(h^{-1})B_{t}^{\rho_{i}}(h).$$

Now we give the explicit formula of Givental quantization of power series of $z$ formulated in the Appendix C of \cite{Ts10}.
\begin{eqnarray*}
&& \bigl(A^i_{t+1}z^m\bigr)^{\wedge} \\
&=& \sum_{\alpha,\beta\in \rep(G)} (E^i_{t+1})^\alpha_\beta  \left( -\sum_{\ell\geq 0}
q^\beta_l \frac{\partial}{\partial q^\alpha_{\ell +t}}
+\frac{\hbar}{2\nu_{\beta}}\sum_{\ell=0}^{t-1} \sum_{\alpha,\beta\in \rep(G)} (-1)^{\ell+1+t}
\frac{\partial}{\partial q^\alpha_\ell}\frac{\partial}{\partial q^\beta_{t-1-\ell} }\right).
\end{eqnarray*}
We have the dilaton shift
$$
\sum_{\beta \in \rep(G)} q_1^\beta\phi_\beta =\sum_{\beta \in G^*}u_1^\beta\phi_\beta -\mathbf{1}
=\sum_{\beta \in \rep(G)}(u_1^\beta-1)\phi_\beta.
$$
So
$$
q_1^\beta = u_1^\beta-1,\quad \forall \beta \in \rep(G),
$$
and
$$
q_k^\beta = u_k^\beta,\quad \forall \beta \in \rep(G), k\neq 1 .
$$
Therefore, we have
\begin{eqnarray*}
\bigl(A^i_{t+1}z^m\bigr)^{\wedge} &=&
\sum_{\alpha,\beta\in \rep(G)} (E^i_{t+1})^\alpha_\beta  \left(\frac{\partial}{\partial u^\alpha_t}-  \sum_{\ell\geq 0}
u^\beta_\ell \frac{\partial}{\partial u^\alpha_{\ell +t}} \right.\\
&& \quad\quad + \left.\frac{\hbar}{2\nu_\beta}\sum_{\ell=0}^{t-1} \sum_{\alpha,\beta\in \rep(G)} (-1)^{\ell+1+t}
\frac{\partial}{\partial u^\alpha_\ell}\frac{\partial}{\partial u^\beta_{t-1-\ell} }\right).
\end{eqnarray*}
\subsection{Sum Over Graphs}
In this section, we write the right hand side of \eqref{eqn:QRR} as a sum
over graphs, following \cite{DSS, DOSS}.

Given a connected graph $\Ga$, we introduce the following notation.
\begin{enumerate}
\item $V(\Ga)$ is the set of vertices in $\Ga$.
\item $E(\Ga)$ is the set of edges in $\Ga$.
\item $H(\Ga)$ is the set of half edges in $\Gamma$.
\item $L^o(\Ga)$ is the set of ordinary leaves in $\Ga$.
\item $L^1(\Ga)$ is the set of dilaton leaves in $\Ga$.
\end{enumerate}

With the above notation, we introduce the following labels:
\begin{enumerate}
\item (genus) $g: V(\Ga)\to \Z_{\geq 0}$.
\item (marking) $\alpha: V(\Ga) \to \rep(G)$. This induces
$\alpha:L(\Ga)=L^o(\Ga)\cup L^1(\Ga)\to \rep(G)$, as follows:
if $l\in L(\Ga)$ is a leaf attached to a vertex $v\in V(\Ga)$,
define $\alpha(l)=\alpha(v)$.
\item (height) $k: H(\Ga)\to \Z_{\geq 0}$.
\end{enumerate}

Given an edge $e$, let $h_1(e),h_2(e)$ be the two half edges associated to $e$. The order of the two half edges does not affect the graph sum formula in this paper. Given a vertex $v\in V(\Ga)$, let $H(v)$ denote the set of half edges
emanating from $v$. The valency of the vertex $v$ is equal to
the size of the set $H(v)$: $\val(v)=|H(v)|$.
A labeled graph $\vGa=(\Ga,g,\alpha,k)$ is {\em stable} if
$$
2g(v)-2 + \val(v) >0
$$
for all $v\in V(\Ga)$.

Let $\bGa(\BG)$ denote the set of all stable labeled graphs
$\vGa=(\Gamma,g,\alpha,k)$. The genus of a stable labeled graph
$\vGa$ is defined to be
$$
g(\vGa):= \sum_{v\in V(\Ga)}g(v)  + |E(\Ga)|- |V(\Ga)|  +1
=\sum_{v\in V(\Ga)} (g(v)-1) + (\sum_{e\in E(\Gamma)} 1) +1.
$$
Define
$$
\bGa_{g,n}(\BG)=\{ \vGa=(\Gamma,g,\alpha,k)\in \bGa(\BG): g(\vGa)=g, |L^o(\Ga)|=n\}.
$$
Define the operator $R(z)$ to be
$$
R(z)=   \exp\Bigl(\sum_{t\geq 1}\frac{(-1)^{t}}{t(t+1)} \sum_{i=1}^m A^i_{t+1} (\frac{z}{\sw_{\rho_{i}}})^t \Bigr).
$$
Then \eqref{eqn:QRR} can be written as
$$
\cD^\tw =R(z)^\wedge  \, \cD^{\BG}.
$$
Under the canonical basis $\{\phi_\gamma:\gamma\in \rep(G)\}$, we have
$$
R(z)^{\alpha}_{\beta}=\frac{|V_{\alpha}|}{|G||V_{\beta}|}\sum_{[h]\in G}\chi_{\alpha}(h^{-1})\chi_{\beta}(h) \exp \Bigl(\sum_{i=1}^{m}\sum_{t=1}^{\infty}\frac{(-1)^t}{t(t+1)} B_{t}^{\rho_i}(h)(\frac{z}{\sw_{\rho_{i}}})^{t}\Bigr).
$$
Given $\beta\in \rep(G)$, define
$$
\bu^\beta(z) = \sum_{a\geq 0} u^\beta_a z^a.
$$

We assign weights to leaves, edges, and vertices of a labeled graph $\vGa\in \bGa(\BG)$ as follows.
\begin{enumerate}
\item {\em Ordinary leaves.}To each ordinary leaf $l \in L^o(\Ga)$ with  $\alpha(l)= \alpha\in \rep(G)$
and  $k(l)= k\in \Z_{\geq 0}$, we assign:
$$(\tilde{\mathcal{L}}^{u})^{\alpha}_{k}(l) = [z^k](\sum_
{\beta}\frac{R(-z)^{\beta}_{\alpha}u^{\beta}(z)}{\sqrt{\nu_{\beta}}}).$$

\item {\em Dilaton leaves.} To each dilaton leaf $l \in L^1(\Ga)$ with $\alpha(l)=\alpha \in \rep(G)$
and $2\leq k(l)=k \in \Z_{\geq 0}$, we assign:

$$(\tilde{\mathcal{L}}^{1})^{\alpha}_{k}(l) = [z^{k-1}](\sum_
{\beta}\frac{R(-z)^{\beta}_{\alpha}}{\sqrt{\nu_{\beta}}}).$$

\item {\em Edges.} To an edge connected a vertex marked by $\alpha\in \rep(G)$ to a vertex
marked by $\beta\in \rep(G)$ and with heights $k$ and $l$ at the corresponding half-edges, we assign:

$$\tilde{\mathcal{E}}^{\alpha,\beta}
_{k,l}(e) = [z^{k}w^{l}]\frac{\delta_{\alpha,\beta}-\sum_{\gamma}R^{\gamma}_{\alpha}(-z)R^{\gamma}_{\beta}(-w)}{z + w}.
$$

\item {\em Vertices.} To a vertex $v$ with genus $g(v)=g\in \Z_{\geq 0}$ and with
marking $\alpha(v)=\gamma\in \rep(G)$, with $n$ ordinary
leaves and half-edges attached to it with heights $k_1, ..., k_n \in \Z_{\geq 0}$ and $m$ more
dilaton leaves with heights $k_{n+1}, \ldots, k_{n+m}\in \Z_{\geq 0}$, we assign

$$\mathcal{V}(v)=\sqrt{v_{\alpha}}^{2-2g(v)-val(v)}\int_{\overline{\mathcal{M}}_{g,n+m}}\psi_{1}^{k_1}\dots\psi_{n+m}^{k_{n+m}}.
$$
\end{enumerate}
We define the weight of a label graph $\overrightarrow{\Gamma}\in\mathbf{\Gamma}(BG)$ to be

$$\widetilde{w}(\overrightarrow{\Gamma})=\prod_{v\in V(\overrightarrow{\Gamma})}V(v)\prod_{e\in E(\overrightarrow{\Gamma})}\tilde{\mathcal{E}}^{\alpha(v_1(e)),\alpha(v_2(e))}
_{k(h_1(e)),k(h_2(e))}(e)\prod_{l\in L^{o}(\overrightarrow{\Gamma})}(\tilde{\mathcal{L}}^{u})^{\alpha(l)}_{k(l)}(l)\prod_{l\in L^{1}(\overrightarrow{\Gamma})}(\tilde{\mathcal{L}}^{1})^{\alpha(l)}_{k(l)}(l)$$
Then
$$
\log(\cD^\tw(\bu))=\log(\hat{R} \cD^{\BG}(\bu)) =
\sum_{\vGa\in \bGa(\BG)}  \frac{ \hbar^{g(\vGa)-1} \widetilde{w}(\vGa)}{|\Aut(\vGa)|}
= \sum_{g\geq 0}\hbar^{g-1} \sum_{n\geq 0}\sum_{\vGa\in \bGa_{g,n}(\cB G)}\frac{\widetilde{w}(\vGa)}{|\Aut(\vGa)|}.
$$
To obtain the graph sum of $\log(\cD^\cX(\hat\bu))$, we consider the fact $\log(\cD^\cX(\hat\bu,\hbar))=\log(\cD^{\tw}(\bu,\be_{1}\hbar))$. We assign a new weight $w(\vGa)$
to each labeled graph $\vGa\in \bGa(\BG)$.  Define
\begin{eqnarray*}
&& (\cL^\bu)^\alpha_k(l) = \frac{1}{\sqrt{\be_1}} (\tcL^\bu)^\alpha_k(l), \quad
(\cL^1)^\alpha_k(l)  = \frac{1}{\sqrt{\be_1}} (\tcL^1)^\alpha_k(l), \\
&& \cE^{\alpha,\beta}_{k,l}(e) = [z^k w^l]
\Bigl(\frac{1}{z+w} (\delta_{\alpha,\beta}-\sum_{\gamma\in G^*} R(-z)^\gamma_\alpha R(-w)^\gamma_\beta)\Bigr),\\
&& \bar{\mathcal{V}}(v)=\sqrt{\bar{v}_{\alpha}}^{2-2g(v)-val(v)}\int_{\overline{\mathcal{M}}_{g,n+m}}\psi_{1}^{k_1}\dots\psi_{n+m}^{k_{n+m}}
\end{eqnarray*}
and define
\begin{eqnarray*}
w(\vGa) &=& \prod_{v\in V(\Ga)} \bar{\mathcal{V}}(v)
\prod_{e\in E(\Ga)}\cE^{\alpha(v_1(e)),\alpha(v_2(e))}_{k(h_1(e)),k(h_2(e))}(e)\\
&& \cdot \prod_{l\in L^o(\Ga)}(\cL^{\bu})^{\alpha(l)}_{k(l)}(l)
\prod_{l\in L^1(\Ga)}(\cL^1)^{\alpha(l)}_{k(l)}(l).
\end{eqnarray*}
With the above notation, we have the following graph sum for $\log(\cD^{\cX}(\hat\bu))$.
\begin{thm}\label{main}
\begin{equation}
\log(\cD^\cX(\hat\bu))
= \sum_{g\geq 0}\hbar^{g-1} \sum_{n\geq 0}\sum_{\vGa\in \bGa_{g,n}(\BG)}\frac{w(\vGa)}{|\Aut(\vGa)|}.
\end{equation}
\begin{equation} \label{eqn:n-point}
\frac{\langle \hat\bu,\ldots,\hat\bu\rangle_{g,n}^{\cX}}{n!} =\sum_{\vGa\in\bGa_{g,n}(\BG)}\frac{w(\vGa)}{|\Aut(\vGa)|}.
\end{equation}
\end{thm}

Finally, we state a straightforward generalization of \eqref{eqn:n-point}. Let
$\bGa_{g,n,n'}(\BG)$ be the set of genus $g$ stable graphs
with $n$ {\em ordered} ordinary leaves and $n'$ {\em unordered} ordinary leaves.
In particular, $\bGa_{g,0,n}(\BG)=\bGa_{g,n}(\BG)$.
Given $\vGa\in \bGa_{g,n,n'}(\BG)$, let $L^O(\Gamma)=\{l_1\ldots, l_n\}$
be the set of ordered ordinary leaves in $\Gamma$, and
let $L^o(\Gamma)$ be the set of unordered ordinary leaves in $\Gamma$;
we assign the following weight to this labeled graph $\vGa\in \Gamma_{g,n,n'}(\BG)$:
\begin{eqnarray*}
w(\vGa) &=& \prod_{v\in V(\Ga)} \bar{\mathcal{V}}(v)
\prod_{e\in E(\Ga)}\cE^{\alpha(v_1(e)),\alpha(v_2(e))}_{k(h_1(e)),k(h_2(e))}(e)\\
&& \cdot \prod_{j=1}^n(\cL^{\bu_j})^{\alpha(l_j)}_{k(l_j)}(l_j)
\prod_{l\in L^o(\Ga)}(\cL^{\bu})^{\alpha(l)}_{k(l)}(l)
\prod_{l\in L^1(\Ga)}(\cL^1)^{\alpha(l)}_{k(l)}(l),
\end{eqnarray*}
where
$$
(\cL^{\bu_j})^\alpha_k(l_j) =  [z^k]\frac{1}{\sqrt{\be_1}}(\sum_
{\beta}\frac{R(-z)^{\beta}_{\alpha}u^{\beta}_{j}(z)}{\sqrt{\nu_{\beta}}}),
\quad
\bu^\beta_j (z) = \sum_{a\geq 0} (u_j)^\beta_a z^a.
$$
Define
$$
\hbu_j =\sum_{\beta\in \rep(G)}\bu^\beta_j(z) \bar{\phi}_\beta.
$$
We have the following generalization of \eqref{eqn:n-point}.
\begin{thm} \label{ordered}
$$
\frac{1}{(n')!}\langle \hbu_1,\ldots,\hbu_n,\bu,\ldots,\bu \rangle_{g,n+n'}^{\cX} =\sum_{\vGa\in\bGa_{g,n,n'}(\BG)}\frac{w(\vGa)}{|\Aut(\vGa)|}.
$$
\end{thm}

\begin{example}[Type D example]
Consider the non-toric Calabi-Yau 3-orbifold $[\C^{3}/D_{n}]$, with $D_{n}$ being binary dihedral group of order $4n$ and $D_{n}$ acting trivially on the third direction. Let $D_{n}=\{<a,b>|a^n=b^2,a^{2n}=1,bab^{-1}=a^{-1}\}$,the action on the $\C^{2}$ component is defined the be:

\begin{equation}
    \rho_{1}(a)=\left(
    \begin{array}{cc}
        \epsilon_{2n} & 0 \\
        0 & \epsilon_{2n}^{-1}
    \end{array}
    \right)
\end{equation}
\begin{equation}
    \rho_{1}(b)=\left(
    \begin{array}{cc}
       0 & i \\
       i & 0
    \end{array}
    \right).
\end{equation}
Here we define $\rho_{1}$ to be the non-trivial 2-dimesional representation and $\rho_{2}$ to be the trivial representation.

So in this case, $A_{t}^{1}$ is the mapping:

$$1_{a^{r}}\to (B_{t}(\frac{r}{2n})+B_{t}(1-\frac{r}{2n}))1_{a^{r}}$$

$$1_{ba^{r}}\to (B_{t}(\frac{1}{4})+B_{t}(\frac{3}{4}))1_{ba^{r}}$$

and $A_{t}^{2}$ is the mapping

$$1_{g}\to B_{t}(0)1_{g}$$

for any $g\in D_{n}$.

The matrix $R(z)^{\alpha}_{\beta}$ is as following:
\begin{eqnarray*}
R(z)^{\alpha}_{\beta}=\frac{|V_{\alpha}|}{|G||V_{\beta}|}\Bigl(\sum_{l=0}^{n}\chi_{\alpha}(a^{-l})\chi_{\beta}(a^{l}) \exp \Bigl(\sum_{t=1}^{\infty}\frac{(-1)^t}{t(t+1)} ((B_{t}(\frac{l}{2n})+B_{t}(1-\frac{l}{2n}))(\frac{z}{\sw_{\rho_{1}}})^{t}+B_{t}(0)(\frac{z}{\sw_{\rho_{2}}})^{t})\Bigr)+\\
(\chi_{\alpha}(b^3)\chi_{\beta}(b)+\chi_{\alpha}(ab^3)\chi_{\beta}(ba))\exp \Bigl(\sum_{t=1}^{\infty}\frac{(-1)^t}{t(t+1)} ((B_{t}(\frac{1}{4})+B_{t}(\frac{3}{4}))(\frac{z}{\sw_{\rho_{1}}})^{t}+B_{t}(0)(\frac{z}{\sw_{\rho_{2}}})^{t})\Bigr)\Bigr).
\end{eqnarray*}

\end{example}
%

\end{document}